\documentclass[twoside,11pt,a4paper]{article}
\usepackage{lineno}
\usepackage{bbm}
\usepackage{bbding}
\usepackage{latexsym}
\usepackage{amsmath}   
\usepackage{tipa}
\usepackage{mathrsfs}  
\usepackage{amssymb}   
\usepackage{amsfonts}
\usepackage{indentfirst}
\usepackage{amsbsy}
\usepackage{amsthm}
\usepackage{bm}
\usepackage{color}
\usepackage{hyperref}
\usepackage{geometry}
\def\beeq{\begin{equation}}      \def\eneq{\end{equation}}
\def\beeqy{\begin{eqnarray}}     \def\eneqy{\end{eqnarray}}
\def\bece{\begin{center}}        \def\ence{\end{center}}
\def\ba {\begin{array}}          \def\ea {\end{array}}
\def\bess{\begin{eqnarray*}}     \def\eess{\end{eqnarray*}}
\def\bes{\begin{split}}          \def\ens{\end{split}}
\def\bali{\begin{align}}         \def\enali{\end{align}}

\def\ali{\aligned}               \def\eali{\endaligned}

\def\pa{\partial} \def\ov{\overline} \def\disp{\displaystyle}
\def\d{\displaystyle\frac}
\def\f{\frac}     \def\q {\quad}      \def\na{\nabla}
     \def\l{\left} \def\r{\right}

\def\s1{\sqrt{-1}}

\def\suml{\sum\limits}

\theoremstyle{plain}
\newtheorem{theorem}{\noindent{\bf Theorem}}[section]
\newtheorem{definition}[theorem]{\noindent{\bf Definition}}
\newtheorem{prop}[theorem]{\noindent{\bf Proposition}}

\newtheorem{cor}[theorem]{\noindent{\bf Corollary}}
\newtheorem{ex}[theorem]{\noindent{\bf Example}}
\newtheorem{rem}[theorem]{\noindent{\bf Remark}}

\def\az{\alpha}      \def\bz{\beta}     \def\gz{\gamma}     \def\dz{\delta}         
\def\zz{\zeta}       \def\ez{\eta}      \def\tz{\theta}              
\def\lz{\lambda}                             
                         
                           \def\oz{\omega}

    \def\BC{{\mathbb{C}}}

    \def\KR{{\mathfrak R}}

          \def\diag {{\rm diag}}

       \def\im {{\rm Im\,}}

   \def\re {{\rm Re\,}}
       \def\tr {{\rm tr\,}}       

        \def\bin#1#2 {{#1\choose#2}}

\def\disp {\displaystyle}
           \def\dfrac#1#2 {{\displaystyle{#1\over#2}}}
           \def\din#1#2 {{\displaystyle{#1\choose#2}}}

\def\disp{\displaystyle}

\let\oldequation\equation
\let\oldendequation\endequation
\renewenvironment{equation}
{\linenomathNonumbers\oldequation}
{\oldendequation\endlinenomath}

\begin{document}
\renewcommand{\thesection}{\arabic{section}}
\renewcommand{\theequation}{\thesection.\arabic{equation}}

\baselineskip=18pt

\title{{\bf The Chern Sectional Curvature of a
Hermitian Manifold}
\author{\small PANDENG CAO, HONGJUN LI$^{*}$\\}
\footnotetext{\scriptsize $^*$Corresponding author.}
}
\date{}
\maketitle

\noindent{\bf Abstract}\hskip3mm
On a Hermitian manifold, the Chern connection can induce a metric connection on the background Riemannian manifold. We call the sectional curvature of the metric connection induced by the Chern connection the Chern sectional curvature of this Hermitian manifold. First, we derive expression of the Chern sectional curvature in local complex coordinates. As an application, we find that a Hermitian metric is K\"ahler if the Riemann sectional curvature and the Chern sectional curvature coincide.  As subsequent results, Ricci curvature and scalar curvature of the metric connection induced by the Chern connection are obtained.
\vspace*{2mm}


\vspace*{2mm}
\noindent{\bf 2020 MR Subject Classification:}\hskip3mm 53C55, 53B35

\thispagestyle{empty}
%
%
%

\section{Introduction}\setcounter{equation}{0}
%
%
%

Suppose $(M,h)$ is an $n$-dimensional Hermitian manifold, and $g=\re h$ is the background Riemannian metric associated to $h$. Let $z=(z^{1},\cdots,z^{n})$ and $x=(x^1,\cdots,x^{2n})$ be local complex and real coordinates of a point $p\in M$, where
$z^\az=x^\az+\sqrt{-1}x^{n+\az}$, $1\leq \az\leq n$.
In this paper, we assume that lowercase Greek indices run from 1 to $n$
and lowercase Latin indices run from 1 to $2n$. In local coordinates, we denote by
$h=h_{\az\bar{\bz}}(z)dz^{\az} d\bar{z}^{\bz}$ and $g=\re{h}=g_{ij}(x){d}x^i{d}x^j$.
Let $\oz=\sqrt{-1}h_{\az\bar{\bz}}(z)dz^{\az}\wedge d\bar{z}^{\bz}$
be the K\"ahler form associated to the Hermitian metric $h$.
If the K\"ahler form $\oz$ is closed, i.e., $d\oz=0$, then we
call $h$ is a K\"ahler metric.

Let us denote by $D$ and $\na$ the Chern connection and the Levi-Civita (or Riemannian) connection, respectively. The curvature operators of the Chern connection and the Levi-Civita connection are denoted by $\KR$ and $R$, respectively.
It is well known that $h$ is K\"ahler if and only if the Chern connection $D$ coincides with the Levi-Civita $\na$ (refer to \cite{Huy,Mo,Zh}, etc.). Hence $\KR$ is the linear extension of $R$ over $\BC$ under K\"ahler hypothesis.

We define a bundle isomorphism $_o:TM\rightarrow T^{1,0}M$ by
\[
 u_o= \frac{1}{2}(u-\sqrt{-1}Ju), \q\forall u\in TM ,
\]
where $J$ is the complex structure on $M$. Let $u,\, v\in TM$,
we set $\xi=u_o,~ \ez=v_o\in T^{1,0}M$. If $h$ is K\"ahler, then
\begin{equation}\label{hbsc}
 R(J{u},u,v,Jv)=2\KR\l(\xi,\bar{\xi},{\ez},\bar{\ez}\r),
 \end{equation}
\begin{equation}\label{rsc}
R(u,v,v,u)=\d{1}{2}\l[\KR\l(\xi,\bar{\ez},{\ez},\bar{\xi}\r)
+\KR\l({\ez},\bar{\xi},\xi,\bar{\ez}\r)
-\KR\l(\xi,\bar{\ez},\xi,\bar{\ez}\r)
-\KR\l(\ez,\bar{\xi},{\ez},\bar{\xi}\r)\r].
\end{equation}
The first formula can be referred to \cite{Zh}, and the second formula
can be referred to \cite{LX,Siu}. Especially, if we take $v=u$ in
\eqref{hbsc} or take $v=J{u}$ in \eqref{rsc}, then \eqref{hbsc}
and \eqref{rsc}
become
\begin{equation}\label{hsc}
R(Ju,u,u,Ju)
 =2\KR(\xi,\bar{\xi},\xi,\bar{\xi})
 \end{equation}
under  K\"ahler hypothesis.

In this paper, we will give a geometric characterization of the right hand
side in \eqref{rsc} when $h$ is non-K\"ahler. In local coordinates,
the right hand side in \eqref{rsc} can be written as
\[
\d{1}{2}\KR_{\az\bar{\bz}\gz\bar{\dz}}
\big(\xi^{\az}\bar{\ez}^{\bz}-{\ez}^{\az}\bar{\xi}^{\bz}\big)
\big({\ez}^{\gz}\bar{\xi}^{\dz}-\xi^{\gz}\bar{\ez}^{\dz}\big),
\]
where $\KR_{\az\bar{\bz}\gz\bar{\dz}}
=\KR\l(\disp\f{\pa}{\pa z^{\az}},\f{\pa}{\pa \bar{z}^{\bz}},
\f{\pa}{\pa z^{\gz}},\f{\pa}{\pa \bar{z}^{\dz}}\r)$.
Under the bundle isomorphism $_o$, any Hermitian connection on $T^{1,0}M$
induces a metric connection on the Riemannian manifold $(M,g)$ \cite{Huy}.
A connection is metric if it compatible with the background Riemannian metric $g=\re h$ \cite{Huy}. 
When no confusion can rise, we still denote by $D$ the metric connection induced by the Chern connection $D$ under the bundle isomorphism $_o$. We find that
\[
\d{\d{1}{2}\KR_{\az\bar{\bz}\gz\bar{\dz}}
\big(\xi^{\az}\bar{\ez}^{\bz}-{\ez}^{\az}\bar{\xi}^{\bz}\big)
\big({\ez}^{\gz}\bar{\xi}^{\dz}-\xi^{\gz}\bar{\ez}^{\dz}\big)}
{\l(h_{\az\bar{\bz}}\xi^{\az}\bar{\xi}^{\bz}\r)\cdot \l(h_{\gz\bar{\dz}}\ez^{\gz}\bar{\ez}^{\dz}\r)
-\d{1}{4}\l[h_{\az\bar{\bz}}\l(\xi^{\az}\bar{\ez}^{\bz}
+\ez^{\az}\bar{\xi}^{\bz}\r)\r]^2}
\]
is just the sectional curvature of the metric connection induced by
the Chern connection $D$ under the bundle isomorphism $_o$, we call it the Chern sectional curvature of $(M,h)$.

\begin{theorem}
Let $(M,h)$ be a Hermitian manifold with the background Riemannian
metric $g=\re h$. Suppose ${D}$ is the metric connection induced by the Chern connection under the bundle isomorphism $_o$.
For arbitrary $u=u^i\f{\pa}{\pa x^i}, ~ v=v^i\f{\pa}{\pa x^i} \in TM$,
we have
\begin{equation}\label{1.4}
g\l(({D}^2u)(v,u),v\r)=\d{1}{2}\KR_{\az\bar{\bz}\gz\bar{\dz}}
\big(\xi^{\az}\bar{\ez}^{\bz}-{\ez}^{\az}\bar{\xi}^{\bz}\big)
\big({\ez}^{\gz}\bar{\xi}^{\dz}-\xi^{\gz}\bar{\ez}^{\dz}\big),
\end{equation}
where $\xi=u_o$ and $\ez=v_o$. Especially,
\begin{equation}\label{1.5}
g\l(({D}^2u)(Ju,u),Ju\r)={2}\KR_{\az\bar{\bz}\gz\bar{\dz}}
\xi^{\az}\bar{\xi}^{\bz}\xi^{\gz}\bar{\xi}^{\dz}.
\end{equation}
\end{theorem}

From \eqref{1.4} in Theorem 1.1, we can see
\[
g\l(({D}^2u)(v,u),v\r)=g\l(({D}^2v)(u,v),u\r).
\]

Recently, Li and Qiu \cite{LQ} proved that a Hermitian
manifold such that \eqref{hsc} is K\"ahler. Hence \eqref{1.5} in
Theorem 1.1 and the main result in \cite{LQ} can imply the following corollary.

\begin{cor}\,
Let $(M,h)$ be a Hermitian manifold such that the Riemann sectional
curvature and the Chern sectional curvature coincide. Then $h$ is
a K\"ahler metric.
\end{cor}

When $h$ is not K\"ahler, there are four Ricci curvatures and two scalar curvatures of the Chern connection $D$, which are respectively denoted by
\[
\ali
\pmb{\rm Ric}^{(1)}_{D} &=\sqrt{-1}\KR^{(1)}_{\az\bar{\bz}}{d}z^{\az}\wedge {d}\bar{z}^{\bz}\q with \q
\KR^{(1)}_{\az\bar{\bz}}=h^{\bar{\dz}\gz}\KR_{\gz\bar{\dz}\az\bar{\bz}},\\
\pmb{\rm Ric}^{(2)}_{D} &=\sqrt{-1}\KR^{(2)}_{\az\bar{\bz}}{d}z^{\az}\wedge {d}\bar{z}^{\bz}\q with \q
\KR^{(2)}_{\az\bar{\bz}}=h^{\bar{\dz}\gz}\KR_{\az\bar{\bz}\gz\bar{\dz}},\\
\pmb{\rm Ric}^{(3)}_{D} &=\sqrt{-1}\KR^{(3)}_{\az\bar{\bz}}{d}z^{\az}\wedge {d}\bar{z}^{\bz}\q with \q
\KR^{(3)}_{\az\bar{\bz}}=h^{\bar{\dz}\gz}\KR_{\gz\bar{\bz}\az\bar{\dz}},\\
\pmb{\rm Ric}^{(4)}_{D} &=\sqrt{-1}\KR^{(4)}_{\az\bar{\bz}}{d}z^{\az}\wedge {d}\bar{z}^{\bz}\q with \q
\KR^{(4)}_{\az\bar{\bz}}=h^{\bar{\dz}\gz}\KR_{\az\bar{\dz}\gz\bar{\bz}},
\eali
\]
\[
\ali
\pmb{s}^{(1)}_D=h^{\bar{\bz}\az}\KR^{(1)}_{\az\bar{\bz}},\q
\pmb{s}^{(2)}_D=h^{\bar{\bz}\az}\KR^{(3)}_{\az\bar{\bz}}.
\eali
\]
We can also define the Ricci curvature of the induced metric connection ${D}$.

Suppose $\{e_i\}_{i=1}^{2n}$ is an orthonormal frame with respect to the background Riemannian metric $g$, we define the Ricci curvature $\pmb{\rm{Ric}}_{{D}}$ of the induced metric connection ${D}$ by
\beeq
\pmb{\rm{Ric}}_{{D}}(u,v)
=\suml_{i=1}^{2n}g\l(({D}^2u)(e_i,v),e_i\r),\q u,\,v \in TM.
\eneq
We denote by
\beeq
{\KR}_{ij}=\suml_{k=1}^{2n}g\l(\l({D}^2\frac{\pa}{\pa x^i}\r)
\l(e_k,\frac{\pa}{\pa x^j}\r),e_k\r)
\eneq
the Ricci curvature tensor of the induced metric connection ${D}$. But we find that
${\KR}_{ij}\neq{\KR}_{ji}$, and
${\KR}_{ij}={\KR}_{ji}$ if and only if the conformal invariant
$b_{\az\bar{\bz}}=0$. The definition of $b_{\az\bar{\bz}}$ can be seen in Section 2. By using Westlake's result \cite{We}, we can see ${\KR}_{ij}={\KR}_{ji}$ if $(M,h)$ is a conformally K\"ahler
manifold.

The scalar curvature $\pmb{s}_{{D}}$ of the induced metric connection ${D}$ is defined by
\beeq
\pmb{s}_{{D}}
=\suml_{i,j=1}^{2n}g\l(({D}^2e_j)(e_i,e_j),e_i\r).
\eneq
We also derive expression of $\pmb{s}_{{D}}$ in local complex
coordinates, and find that $\pmb{s}_{{D}}$ is equal to four times of the second Chern scalar curvature $\pmb{s}^{(2)}_D$.

\section{The Chern sectional curvature}\setcounter{equation}{0}

The Chern connection $D$ is the unique connection on the holomorphic
tangent bundle $T^{1,0}M$ which is compatible with the Hermitian metric
$h$ and the complex structure $J$.
We denote by
$$
\tz=(\tz_{\bz}^{\az})=\l({\varGamma}_{\bz;\gz}^{\az}dz^{\gz}\r)
=\l(
h^{\bar{\lz}\az}\d{\pa h_{\bz\bar{\lz}}}{\pa z^{\gz}}dz^\gz\r)
$$
the $n \times n$ matrix of connection 1-forms.
The curvature operator $\KR$ of $D$ is defined by
\begin{equation}
\KR(\xi,\bar{\ez},\zz,\bar{\chi})
=h\left(\left(D_\zz D_{\bar\chi}-D_{\bar\chi}D_\zz
-D_{[\zz,\bar{\chi}]}\right)\xi,\bar{\ez}\right), \q
\forall \xi,\;\ez,\;\zz,\;\chi \in T^{1,0}M,
\end{equation}
where $[\,,\,]$ is the Lie bracket.
We call
\begin{equation}\label{2.1}
\KR_{\az\bar{\bz}\gz\bar{\dz}}
=\KR\l(\disp\f{\pa}{\pa z^{\az}},\f{\pa}{\pa \bar{z}^{\bz}},
\f{\pa}{\pa z^{\gz}},\f{\pa}{\pa \bar{z}^{\dz}}\r)
=-\d{\pa^2 h_{\az\bar{\bz}}}{\pa z^{\gz}\pa \bar{z}^{\dz}}+
\d{\pa h_{\az\bar{\lz}}}{\pa z^{\gz}}h^{\bar{\lz}\kappa}
\d{\pa h_{\kappa\bar{\bz}}}{\pa \bar{z}^{\dz}}
\end{equation}
the holomorphic sectional curvature tensor of $D$.
Ulteriorly, the holomorphic sectional curvature is defined by
\begin{equation}\label{2.3}
{\rm HSC}(\xi)
=\d{\KR\l(\xi,\bar{\xi},\xi,\bar{\xi}\r)}
{h\l(\xi,{\xi}\right)^2},\q \forall\xi\in T^{1,0}M.
\end{equation}

We still denote by ${D}$ the metric connection on the Riemannian manifold $(M,g)$, which is induced by the Chern connection $D$ under
the bundle isomorphism $_o$. Then
\begin{equation}\label{2.4}
\ali
&{D}\f{\pa}{\pa x^{\az}}
=\f{1}{2}\l(\tz_\az^\bz+\bar{\tz}_\az^\bz\r)\otimes\f{\pa}{\pa x^{\bz}}
-\f{\sqrt{-1}}{2}\l(\tz_\az^\bz-\bar{\tz}_\az^\bz\r)\otimes\f{\pa}{\pa x^{\bz+n}},
\\
&{D}\f{\pa}{\pa x^{\az+n}}
=\f{\sqrt{-1}}{2}\l(\tz_\az^\bz-\bar{\tz}_\az^\bz\r)
\otimes\f{\pa}{\pa x^{\bz}}
+\f{1}{2}\l(\tz_\az^\bz+\bar{\tz}_\az^\bz\r)
\otimes\f{\pa}{\pa x^{\bz+n}}.
\eali
\end{equation}
For a general Hermitian manifold $(M,h)$, the induced metric connection
${D}$ is compatible with the background Riemannian metric $g$, but not torsion free.

\begin{definition}
Let $(M,h)$ be a Hermitian manifold.
We call
\begin{equation}
\begin{split}
K_{{D}}(u,v)
=\d{g(({D}^2u)(v,u),v)}{g(u,u)g(v,v)-g(u,v)^2}
\end{split}
\end{equation}
the Chern sectional curvature of the 2-plane $\Pi(u,v)$ spanned by
two linearly independent tangent vectors $u=u^i\f{\pa }{\pa x^i},\,
v=v^i\f{\pa }{\pa x^i}\in TM$.
\end{definition}

The Levi-Civita connection $\nabla$ is the unique connection on the tangent bundle $TM$, which is torsion free and compatible with the background Riemannian metric $g=\re h$. Let us denote by
\[
\varphi_{i}^{k}=\gz^{k}_{ij}dx^j
\]
connection 1-forms of the Levi-Civita connection $\nabla$,
where
\[
\gz_{ij}^{k}=\d{1}{2}g^{kl}\l(\d{\pa g_{il}}{\pa x^{j}}+\d{\pa g_{jl}}{\pa x^{i}}-\d{\pa g_{ij}}{\pa x^{l}}\r)
\]
is the connection coefficient of $\nabla$.
The curvature operator $R$ of $\nabla$ is defined by
\beeq
R(u,v,w,y)=g\l(\l(\na_u\na_v-\na_v\na_u-\na_{\l[u,v\r]}\r)w,y\r),\q
\forall u,\; v,\;w,\; y\in TM.
\eneq
The Riemann sectional curvature on the 2-plane $\Pi(u,v)$ is defined by
\begin{equation}
K_{\na}(u,v)
=\d{R(u,v,v,u)}{g(u,u)g(v,v)-g(u,v)^2}.
\end{equation}
We call $\Pi(u, Ju)$ a holomorphic plane section \cite{Zh}.

If $(M,h)$ is a K\"ahler manifold, then the induced metric connection ${D}$ and the Levi-Civita connection $\nabla$ coincide, thus $K_{{D}}(u,v)=K_{\na}(u,v)$.

Now we extend the background Riemannian metric $g$ linearly over $\BC$ to the complexified tangent bundle $T_{\BC}M=TM\otimes\BC=T^{1,0}M\oplus T^{0,1}M$, then
$$
\ali
g\l(\frac{\pa}{\pa z^{\az}}, \frac{\pa}{\pa z^{\bz}}\r)
&=g\l(\frac{\pa}{\pa \bar{z}^{\az}}, \frac{\pa}{\pa \bar{z}^{\bz}}\r)=0,\\
g\l(\frac{\pa}{\pa z^{\az}}, \frac{\pa}{\pa \bar{z}^{\bz}}\r)
&=\ov{g\l(\frac{\pa}{\pa \bar{z}^{\az}}, \frac{\pa}{\pa {z}^{\bz}}\r)}=\frac{1}{2}h_{\az\bar{\bz}}.
\eali
$$
Hence
\beeq
h\l(\xi,\bar\ez\r)=2g\l(\xi,\bar\ez\r),\q\forall \xi,\,\ez\in T^{1,0}M.
 \eneq

Next we give expression of the Chern sectional curvature in local
complex coordinates.

{\bf Proof of Theorem 1.1}\q
We denote by ${D}\f{\pa}{\pa x^i}=\tilde{\tz}_i^j\f{\pa}{\pa x^j}$
and $\tilde{\tz}=(\tilde{\tz}_i^j)$. It follows from \eqref{2.4} that
\begin{equation}
\tilde{\tz}=F\diag\l\{\tz,\bar\tz\r\}F^{-1},
\end{equation}
 where $F = \left(
\begin{array}{cr}
I&I\\
\sqrt {-1}I& -\sqrt{-1}I
\end{array} \right)$, $F^{-1}=\displaystyle\frac{1}{2}\left(
\begin{array}{cr}
I&-\sqrt{-1}I\\
I& \sqrt{-1}I
\end{array} \right)$ is the inverse of $F$, $I$ is the $ n\times n$ identity matrix, $\diag\l\{\tz,\bar{\tz}\r\}=\left(\begin{array}{cr}
\tz&0\\
0& \bar{\tz}
\end{array} \right)$.
In order to simplify the calculation process, we introduce the
following notations. Set
$$
\f{\pa}{\pa x}=\l(\frac{\pa}{\pa x^1},\cdots,\frac{\pa}{\pa x^{2n}}\r),\q
\f{\pa}{\pa z}=\l(\frac{\pa}{\pa z^1},\cdots,\frac{\pa}{\pa z^{n}}\r),\q
\f{\pa}{\pa \bar{z}}=\l(\frac{\pa}{\pa \bar{z}^1},\cdots,
\frac{\pa}{\pa \bar{z}^{n}}\r),
$$
then
\[
\f{\pa}{\pa x}=\left(\f{\pa}{\pa z}, \f{\pa}{\pa \bar{z}}\right)F^{\rm t},
\]
where $F^{\rm t}$ means the transpose of $F$. We still denote by $A^{\rm t}$ the transpose of a vector $A$.
For any $u=u^i\frac{\pa}{\pa x^i}\in TM$, $\xi=u_o=\xi^\az\frac{\pa }{\pa z^\az}\in T^{1,0}M$, we set
$$
\pmb{u}=(u^1,\cdots,u^{2n}),\q
\pmb{\xi}=(\xi^1,\cdots,\xi^n),
$$
then
$$
\pmb{u}=(\pmb{\xi},\bar{\pmb{\xi}})F^{-1}.
$$
By a straightforward computation, we have
\[
\ali
{D}^2u=&u^j\l(d\tilde{\tz}_j^i-\tilde{\tz}_j^k\wedge\tilde{\tz}_k^i\r)
\otimes\f{\pa}{\pa x^i}\\
=&\pmb{u}\l(d\tilde{\tz}-\tilde{\tz}\wedge\tilde{\tz}\r)
\otimes\l(\f{\pa}{\pa x}\r)^{\rm t}\\
=&(\pmb{\xi},\bar{\pmb{\xi}})\diag\l\{d\tz-\tz\wedge \tz,
d\bar\tz-\bar\tz\wedge \bar\tz\r\}\otimes
\left(
 \begin{array}{c}
\l(\f{\pa }{\pa z}\r)^{\rm t} \\
\l(\f{\pa }{\pa \bar{z}}\r)^{\rm t} \\
\end{array}
 \right)\\
=&\pmb{\xi}\l(\bar{\pa}\tz\r)\otimes\l(\f{\pa }{\pa z}\r)^{\rm t}
+\bar{\pmb{\xi}}\l({\pa}\bar{\tz}\r)
\otimes\l(\f{\pa }{\pa \bar{z}}\r)^{\rm t}\\
=&\xi^\az\l(\bar{\pa}\tz_\az^\mu\r)\otimes\f{\pa }{\pa z^{\mu}}
+\bar{\xi}^{\bz}\l({\pa}\bar{\tz}_\bz^\mu\r)
\otimes\f{\pa }{\pa \bar{z}^{\mu}}\\
=&\xi^\az\l(\KR^{\mu}_{\az\gz\bar{\dz}}dz^\gz\wedge d\bar{z}^\dz\r)\otimes
\f{\pa}{\pa z^\mu}
-\bar{\xi}^\bz\l(\ov{\KR^{\mu}_{\bz\dz\bar{\gz}}}dz^\gz\wedge
d\bar{z}^\dz\r)\otimes\f{\pa}{\pa \bar{z}^\mu},
\eali
\]
where $d=\pa+\bar{\pa}$ is the exterior differentiation.
Hence
\[
\big({D}^2u\big)(v,w)=\big({\ez}^\gz\bar{\zz}^\dz
-\zz^\gz \bar{\ez}^{\dz}\big)\l(\xi^\az\KR^{\mu}_{\az\gz\bar{\dz}}
\f{\pa}{\pa z^\mu}
-\bar{\xi}^\bz\ov{\KR^{\mu}_{\bz\dz\bar{\gz}}}\f{\pa}{\pa \bar{z}^\mu}\r),
\]
\beeq\label{2.9}
\ali
g\l(\big({D}^2u\big)(v,w),y\r)
=&\f{1}{2}\KR_{\az\bar{\bz}\gz\bar{\dz}}\big(\xi^\az\bar{\chi}^{\bz}
-\chi^\az\bar{\xi}^{\bz}\big)\big({\ez}^\gz\bar{\zz}^\dz-\zz^\gz \bar{\ez}^{\dz}\big),
\eali
\eneq
where $w,\, y\in TM$, $\zeta=w_o, ~\chi=y_o \in T^{1,0}M$.
Especially, we can obtain \eqref{1.4} by taking $w=u$ and $y=v$ in \eqref{2.9}. Note that
\[\l(Ju\r)_o=\d{1}{2}\l(Ju-\sqrt{-1}J^2u\r)
=\d{\sqrt{-1}}{2}\l(u-\sqrt{-1}J u\r)=J(u_o).
\]
In order to prove \eqref{1.5}, we only replace $v$ and $\ez$ in \eqref{1.4} with $Ju$ and $\sqrt{-1}v$, respectively.\qed

According to \eqref{1.4} in Theorem 1.1, we have
\beeq
K_{{D}}(u,v)=\d{\d{1}{2}\KR_{\az\bar{\bz}\gz\bar{\dz}}
\big(\xi^{\az}\bar{\ez}^{\bz}-{\ez}^{\az}\bar{\xi}^{\bz}\big)
\big({\ez}^{\gz}\bar{\xi}^{\dz}-\xi^{\gz}\bar{\ez}^{\dz}\big)}
{\l(h_{\az\bar{\bz}}\xi^{\az}\bar{\xi}^{\bz}\r)\cdot \l(h_{\gz\bar{\dz}}\ez^{\gz}\bar{\ez}^{\dz}\r)
-\d{1}{4}\l[h_{\az\bar{\bz}}\l(\xi^{\az}\bar{\ez}^{\bz}
+\ez^{\az}\bar{\xi}^{\bz}\r)\r]^2}.
\eneq
It follows from the above formula that
\[
K_{{D}}(u,v)=K_{{D}}(v,u).
\]

\begin{rem}
By \eqref{2.9}, we have
\beeq
g\l(\big({D}^2u\big)(Jw,w),Ju\r)
=2\KR\l(u_o,\ov{u_o},w_o,\ov{w_o}\r).
\eneq
\end{rem}

\begin{rem}
By \eqref{2.9}, it is clear that
\[
g\l(\big({D}^2u\big)(v,w),y\r)
=g\l(\big({D}^2u\big)(Jv,Jw),y\r)
=g\l(\big({D}^2Ju\big)(v,w),Jy\r),
\]
\[
g\l(\big({D}^2u\big)(v,w),y\r)
=-g\l(\big({D}^2y\big)(v,w),u\r)
=-g\l(\big({D}^2u\big)(w,v),y\r).
\]
But we can not expect the following formula
\[
g\l(\big({D}^2u\big)(v,w),y\r)
=g\l(\big({D}^2v\big)(u,y),w\r)
\]
always holds for arbitrary $u,\, v, \, w, \, y\in TM$.
\end{rem}

As an application, we can extend Lu's result \cite{Lu} of K\"ahler manifolds with non-negative (or non-positive) Riemann sectional curvatures on K\"ahler-like manifolds with non-negative (or non-positive) Chern sectional curvature. Yang and Zheng \cite{YZ} defined K\"ahler-like manifolds, which are classes of non-K\"ahler manifolds.

\begin{definition}{\rm\cite{YZ}}
A Hermitian metric $h$ is called K\"ahler-like, if its holomorphic sectional curvature tensor satisfies $\mathfrak{R}_{\az \bar{\bz}\gz\bar{\dz}}
=\mathfrak{R}_{\gz \bar{\bz} \az\bar{\dz}}$ for all
$\az,\, \bz, \, \gz, \, \dz=1,\,2,\cdots, n$.
\end{definition}

By using the same method as that in \cite{Lu}, we have the following result without details.

\begin{prop}
Let $(M,h)$ be a K\"ahler-like manifold with non-negative (resp.
non-positive) Chern sectional curvature. Then
\begin{equation}
\l|\KR\l(\xi,\bar{\xi},\ez,\bar{\ez}\r)\r|^2\leq
\KR\l(\xi,\bar{\xi},\xi,\bar{\xi}\r)\KR\l(\ez,\bar{\ez},\ez,\bar{\ez}\r).
\end{equation}
\end{prop}

\begin{prop}
For arbitrary $u=u^i\f{\pa}{\pa x^i},~w=w^i\f{\pa}{\pa x^i}\in TM$,
we have
\begin{equation}\label{2.13}
\pmb{\rm{Ric}}_{{D}}(u,w)
=\KR^{(4)}_{\az\bar{\dz}}\xi^{\az}\bar{\zz}^{\dz}
+\KR^{(3)}_{\az\bar{\dz}}\zz^{\az}\bar{\xi}^{\dz},
\end{equation}
where $\xi=u_o$, $\zz=w_o$.
\end{prop}
\begin{proof}
 For any lowercase Greek index $\az$, we denote by $\az^*=\az+n$.
By \eqref{2.9}, we have
\[
\ali
g\l(\big({D}^2u\big)\l(\frac{\pa}{\pa x^{\kappa}},w\r),\frac{\pa}{\pa x^{\lz}}\r)
=&\re\big(\KR_{\az\bar{\lz}\kappa\bar{\dz}}\xi^\az\bar{\zz}^\dz
-\KR_{\az\bar{\lz}\gz\bar{\kappa}}\xi^\az{\zz}^\gz\big),\\
\eali
\]
\[
\ali
g\l(\big({D}^2u\big)\l(\frac{\pa}{\pa x^{\kappa}},w\r),\frac{\pa}{\pa x^{\lz^*}}\r)
=&\im\big(\KR_{\az\bar{\lz}\kappa\bar{\dz}}\xi^\az\bar{\zz}^\dz
- \KR_{\az\bar{\lz}\gz\bar{\kappa}}\xi^\az{\zz}^\gz
\big),\\
\eali
\]
\[
\ali
g\l(\big({D}^2u\big)\l(\frac{\pa}{\pa x^{\kappa^*}},w\r),\frac{\pa}{\pa x^{\lz}}\r)
=&-\im\big(\KR_{\az\bar{\lz}\kappa\bar{\dz}}\xi^\az\bar{\zz}^\dz
+\KR_{\az\bar{\lz}\gz\bar{\kappa}}\xi^\az{\zz}^\gz\big),\\
\eali
\]
\[
\ali
g\l(\big({D}^2u\big)\l(\frac{\pa}{\pa x^{\kappa^*}},w\r),
\frac{\pa}{\pa x^{\lz^*}}\r)
=&\re\big(\KR_{\az\bar{\lz}\kappa\bar{\dz}}\xi^\az\bar{\zz}^\dz
+\KR_{\az\bar{\lz}\gz\bar{\kappa}}\xi^\az{\zz}^\gz\big).\\
\eali
\]

Set $L=\l(L_{\kappa\bar{\lz}}\r)_{1\leq\kappa,\lz\leq n}=L_1+\sqrt{-1}L_2$,
where $L_{\kappa\bar{\lz}}=\KR_{\az\bar{\lz}\kappa\bar{\dz}}
\xi^\az\bar{\zz}^\dz$, $L_1=\re L$ and $L_2=\im L$.
Set $K=\l(K_{\kappa {\lz}}\r)_{1\leq\kappa,\lz\leq n}=K_1+\sqrt{-1}K_2$,
where $L_{\kappa{\lz}}=\KR_{\az\bar{\lz}\gz\bar{\kappa}}
\xi^\az \zz^\gz$, $K_1=\re K$ and $K_2=\im K$.
Being similar to the formula (2.18) in \cite{LQZ}, we can write as
\[
\ali
&\left(
              \begin{array}{cc}
                L_1-K_1 & L_2-K_2\\
                -L_2-K_2 & L_1+K_1 \\
              \end{array}
            \right)\\
            =&F\diag\l\{L ,\bar{L} \r\}F^{-1}
-\frac{1}{2}\diag\l\{I,-I\r\}F\diag\l\{K,\bar{K}\r\}
F^{\mathrm{t}}\diag\l\{I,-I\r\}.
\eali
\]
Denote by $G=(g_{ij})$ and $H=(h_{\az\bar{\bz}})$, then
\[
G^{-1}=F\diag\l\{H^{-1},\bar{H}^{-1}\r\}F^{-1}.
\]
A direct computation shows that
\[
\ali
\tr\l[G^{-1}\left(
              \begin{array}{cc}
                L_1-K_1 & L_2-K_2\\
                -L_2-K_2 & L_1+K_1 \\
              \end{array}
            \right)\r]
=\tr\l(H^{-1}L\r)+\tr\l(\bar{H}^{-1}\bar{L}\r)
=h^{\bar{\lz}\kappa} \KR_{\az\bar{\lz}\kappa\bar{\dz}} \xi^{\az}\bar{\zz}^{\dz}+
\ov{h^{\bar{\lz}\kappa} \KR_{\az\bar{\lz}\kappa\bar{\dz}} \xi^{\az}\bar{\zz}^{\dz}},
\eali
\]
 where $\tr(\,)$ means the trace of a square matrix.
Hence
\[
\pmb{\rm{Ric}}_{ D}(u,w)
=g^{kl}g\l(\l({D}^2u\r)\l(\frac{\pa}{\pa x^{k}},w\r),
\frac{\pa}{\pa x^{l}}\r)
=\KR^{(4)}_{\az\bar{\dz}}\xi^{\az}\bar{\zz}^{\dz}
+\KR^{(3)}_{\az\bar{\dz}}\zz^{\az}\bar{\xi}^{\dz}.
\]
This completes the proof.
\end{proof}

We recall that the Ricci curvature tensor ${\KR}_{ij}$ of the connection ${D}$ is defined by

\[
{\KR}_{ij}=\suml_{k=1}^{2n}g\l(\l({D}^2\frac{\pa}{\pa x^i}\r)
\l(e_k,\frac{\pa}{\pa x^j}\r),e_k\r)
=g^{kl}g\l(\l({D}^2\frac{\pa}{\pa x^i}\r)
\l(\frac{\pa}{\pa x^{k}},\frac{\pa}{\pa x^j}\r),\frac{\pa}{\pa x^{l}}\r).
\]
From \eqref{2.13}, we can see ${\KR}_{ij}\neq {\KR}_{ji}$ in general.
It follows from \eqref{2.13} that
${\KR}_{ij}={\KR}_{ji}$ if and only if $
\KR^{(3)}_{\az\bar{\dz}}=\KR^{(4)}_{\az\bar{\dz}}$.
Are there some Hermitian manifolds such that ${\KR}_{ij}= {\KR}_{ji}$?

Let
\[
T_{\bz \gz}^{\az}=\varGamma_{\bz;\gz}^{\az}-\varGamma_{\gz;\bz}^{\az}
=h^{\bar{\lz}\az}\l(\d{\pa h_{\bz\bar{\lz}}}{\pa z^{\gz}}
-\d{\pa h_{\gz\bar{\lz}}}{\pa z^{\bz}}\r)
\]
be the torsion tensor of the Chern connection $D$. Set $T_{\bz}=\suml_{\az=1}^nT^\az_{\bz\az}$.
 Lee \cite{Lee} first introduced the following conformal invariant
\beeq
b_{\az\bar{\bz}}=\frac{\pa T_{\az}}{\pa \bar{z}^\bz}
-\frac{\pa T_{\bar{\bz}}}{\pa {z}^\az},
\eneq
where $T_{\bar{\bz}}=\ov{T_{{\bz}}}$.  In fact, it is easy to see
\beeq
b_{\az\bar\bz}=\KR^{(4)}_{\az\bar{\bz}}-\KR^{(3)}_{\az\bar{\bz}}.
\eneq
If there exists a positive scalar function $\rho$ in $(M,h)$ such that $\rho h$ is a K\"ahler metric, we call $(M,h)$ a conformally K\"ahler manifold.
Westlake \cite{We} proved $b_{\az\bar{\bz}}=0$ on conformally K\"ahler manifolds.  Hence ${\KR}_{ij}= {\KR}_{ji}$ on conformally K\"ahler manifolds \cite{Gau84}.
Next we provide an example such that ${\KR}_{ij}= {\KR}_{ji}$.
 \begin{ex}
Let $\mathbb{B}^n=\{z\in \BC^n: |z|^2=\suml_{\az=1}^n|z^\az|^2<1\}$ be the unit ball in $\BC^n\,(n>1)$ endowed with a Hermitian metric $h=\disp\frac{4|dz|^2}{(1-|z|^2)^2}$.
The symbols put as above.
A direct computation shows that
\[
\ali
{\KR}_{\az\bar{\bz}\mu\bar{\nu}}=&\frac{4\dz_{\az\bz}}{(1-|z|^2)^4}\left[
(1-|z|^2)\dz_{\mu\nu}+\bar{z}^{\mu}z^{\nu}\right],\\
\KR^{(3)}_{\az\bar{\bz}}=&\frac{2}{(1-|z|^2)^2}\left[
(1-|z|^2)\dz_{\az\bz}+\bar{z}^{\az}z^{\bz}\right],\\
\KR^{(4)}_{\az\bar{\bz}}=&\frac{2}{(1-|z|^2)^2}\left[
(1-|z|^2)\dz_{\az\bz}+\bar{z}^{\az}z^{\bz}\right],
\eali
\]
where $\dz_{\az\bz}
$ are Kronecker symbols. Hence $\KR^{(4)}_{\az\bar{\bz}}=\KR^{(3)}_{\az\bar{\bz}}$ yields
${\KR}_{ij}= {\KR}_{ji}$.
\end{ex}

\begin{prop}
Let $(M,h)$ be a Hermitian manifold with the background Riemannian
metric $g=\re h$. Suppose ${D}$ is the metric connection induced by the Chern connection under the bundle isomorphism $_o$.
Then $\pmb{s}_{{D}}$ is equal to four times of the second Chern scalar curvature $\pmb{s}^{(2)}_D$, i.e.,
\begin{equation}\label{ }
\pmb{s}_{{D}}=4\pmb{s}^{(2)}_D.
\end{equation}
\end{prop}
\begin{proof}
For simplicity, we denote by
\[
\pmb{\rm{Ric}}_{{D}}(u):=\pmb{\rm{Ric}}_{{D}}(u,u)
={\KR}_{ij}u^iu^j,
\]
then
\[
\ali
\pmb{s}_{{D}}&=\frac{1}{2}g^{ij}\l(
{\KR}_{ij}+{\KR}_{ji}\r)
=\frac{1}{2}g^{ij}\f{\pa^2\pmb{\rm{Ric}}_{{D}}(u)}{\pa u^i \pa u^j}
=2h^{\bar{\dz}\az}\f{\pa^2\pmb{\rm{Ric}}_{{D}}(u)}{\pa \xi^\az \pa \bar{\xi}^\dz}\\
&=2h^{\bar{\dz}\az}h^{\bar{\lz}\kappa}\l(\KR_{\az\bar{\lz}\kappa\bar{\dz}}+
\KR_{\kappa\bar{\dz}\az\bar{\lz}}\r)=4\pmb{s}^{(2)}_D.
\eali
\]
This completes the proof.
\end{proof}

There are various Hermitian connections on the holomorphic tangent bundle in \cite{Gau97}. 
We can consider the metric connection induced by any Hermitian connection under 
the bundle isomorphism $_o$.
For example, the Levi-Civita connection $\na^{lc}$ on $T^{1,0}M$ (the restriction of the complexified Levi-Civita connection $\na$ to $T^{1,0}M$, see \cite{HLY}). 
We denote by $\pmb{s}_{\na^{lc}}$ the scalar curvature of the metric connection induced by the Levi-Civita connection $\na^{lc}$ on $T^{1,0}M$. 
Gauduchon \cite{Gau84} shown that $\pmb{s}_{\na^{lc}}$ and the Riemannian scalar curvature $\pmb{s}_{\na}$ coincide on a compact Hermitian manifold if and only if the compact Hermitian manifold is balanced.

\bigskip

\noindent
{\bf Acknowledge}: This work is supported by the National Natural Science Foundation of China (Grant No. 12001165). The authors are very grateful to the referee for providing many helpful suggestions.

\bigskip

PANDENG CAO

SCHOOL OF MATHEMATICAL SCIENCES

XIAMEN UNIVERSITY

XIAMEN, FUJIAN PROVINCE, 361005, CHINA

Email address: caopandeng2023@163.com

\

HONGJUN LI

SCHOOL OF MATHEMATICAL AND STATISTICS

HENAN UNIVERSITY

KAIFENG, HENAN  PROVINCE, 475004, CHINA

Email address: lihj@vip.henu.edu.cn


\smallskip


\begin{thebibliography}{99}
\bibitem{Gau84}
P. Gauduchon, \emph{La 1-forme de torsion d'une vari$\acute{e}$t$\acute{e}$ hermitienne compacte}, Math. Ann. {\bf}267 (1984), no. 4, 495--518. https://doi.org/10.1007/BF01455968

\bibitem{Gau97}
P. Gauduchon, \emph{Hermitian connections and Dirac operators}.
Boll. Un. Mat. Ital. B (7), {\bf11} (1997), no. 2, suppl., 257--288.

\bibitem{HLY}
J. He,  K. Liu, and X. Yang, \emph{Levi-Civita Ricci-flat metrics on compact complex manifolds}, J. Geom. Anal. {\bf30}(2020), no. 1, 646--666.
https://doi.org/10.1007/s12220-019-00156-9

\bibitem{Huy}
D. Huybrechts, \emph{Complex Geometry: An Introduction}, Springer, Berlin
Heidelber, 2005. https://doi.org/10.1007/b137952

\bibitem{Lee}
H. Lee, \emph{A kind of even dimensional differential geometry and its application to exterior calculus},
Amer. J.  Math. {\bf65} (1943), no. 3, 433--438. https://doi.org/10.2307/2371967

\bibitem{LQ}
H. Li and C. Qiu, \emph{A note on holomorphic sectional curvature of a Hermitian manifold}, Glasgow Math. J. {\bf 64} (2022), no. 3, 739--745.
https://doi.org/10.1017/S0017089522000064

\bibitem{LQZ}
H. Li, C. Qiu, and G. Zhong, \emph{Curvatures of strongly convex K\"ahler-Finsler manifolds}, Differ. Geom. Appl. {\bf 86} (2023), 101957(25pages). https://doi.org/10.1016/j.difgeo.2022.101957

\bibitem{Lu}
Q. Lu, \emph{A note on the extremum of sectional curvature of a K\"ahler manifold}, Sci. Sin. Ser. A  {\bf27} (1984), no. 4, 367--371.

\bibitem{LX}
Q. Lu (K. Look) and Y. Xu(I. Hs\"u),  \emph{A note on transitive domains},
Acta Math. Sin. {\bf11} (1961), no. 1, 11--23; English transl.,
Chin. J. Math. {\bf2} (1962), 11--26. 

\bibitem{Mo}
A. Moroianu, \emph{Lectures on K\"ahler geometry}, Cambridge University Press,
Cambridge, 2007. https://doi.org/10.1017/CBO9780511618666

\bibitem{Siu}
Y. Siu, \emph{The complex-analyticity of harmonic maps and the strong rigidity
 of compact K\"ahler manifolds},  Ann. Math. {\bf112} (1980), no. 1, 73--111.
https://doi.org/10.2307/1971321

\bibitem{We}
W.-J. Westlake, \emph{Conformally K\"ahler manifolds}, Math. Proc. Camb. Phil. Soc. {\bf50} (1954), 16--19. https://doi.org/10.1017/S0305004100029029

\bibitem{YZ}
B. Yang and F. Zheng, \emph{On curvature tensors of Hermitian manifolds},
Commun. Anal. Geom. {\bf26} (2018), no. 5, 1195--1222. https://doi.org/10.4310/CAG.2018.v26.n5.a7

\bibitem{Zh}
F. Zheng, \emph{Complex differential geometry},
American Mathematical Society, International Press, Boston, 2000.
https://doi.org/10.1090/amsip/018


\end{thebibliography}
\end{document}